\documentclass[11pt]{amsart}

\newcommand{\zz}{{\mathbb Z}}
\newcommand{\zc}{{\mathbb C}}
\newcommand{\glnzt}{GL_{n-1}(\zz[t,t^{-1}])}
\newcommand{\glnzq}{GL_{\binom{n}{2}}(\zz[t^{\pm 1},q^{\pm 1}])}
\newcommand{\glnct}{GL_{n-1}(\zc[t,t^{-1}])}
\newcommand{\glncq}{GL_{\binom{n}{2}}(\zc[t^{\pm 1},q^{\pm 1}])}
\newcommand{\ra}{\rightarrow}
\newcommand{\lra}{\longrightarrow}
\newcommand{\lslnc}{{\mathfrak sl}_{n-1}(\zc)}
\newcommand{\galpha}{\Gamma_{n(n-2)}(\alpha)}

\title[Representations of Braid Groups]{Representations of Braid Groups via Conjugation
Actions on Congruence Subgroups} 
  
  
\author{Kevin P.~Knudson} 
\address{Department of Mathematics and Statistics\\
Mississippi State University\\
P.O.~Drawer MA\\
Mississippi State, MS 39762} 
\email{knudson@math.msstate.edu} 
\urladdr{http://www2.msstate.edu/$\sim$kk116/} 
  
\thanks{Partially supported by NSF grant DMS-0242906, and by ORAU} 
  
\newtheorem{theorem}{Theorem}[section] 
\newtheorem{prop}[theorem]{Proposition} 
\newtheorem{lemma}[theorem]{Lemma} 
\newtheorem{cor}[theorem]{Corollary} 
  
\theoremstyle{definition}

\theoremstyle{remark} 
  
\newtheorem{remark}[theorem]{Remark}

\begin{document} 
  
\begin{abstract} 
We construct two families of representations of the braid group $B_n$ by considering
conjugation actions on congruence subgroups of $GL_{n-1}(\zz[t^{\pm 1},q^{\pm 1}])$.  Many of
these representations are shown to be faithful.
\end{abstract} 
  
\maketitle 
  
\section{Introduction}\label{intro} 

Denote by $B_n$ the braid group on $n$ strings.  The purpose of this note is to construct
two families of representations
$$\rho_n(\alpha):B_n\lra SL_{n(n-2)}(\zc)$$
and 
$$\mu_n(\alpha,\beta):B_n\lra SL_N(\zc), \quad N=\binom{n}{2}^2 - 1,$$
where $\alpha$ and $\beta$ are nonzero complex numbers.  If $\alpha,\beta$
lie in some subfield $F$ of $\zc$, then the representations are defined over
$F$ (in fact, over $\zz[\alpha^{\pm 1},\beta^{\pm 1}]$).

The starting point for $\rho_n(\alpha)$ is the reduced Burau representation
$\beta_n:B_n\ra \glnzt$.  For each $i\ge 1$, set
$$K^i(\alpha) = \{A\in SL_{n-1}(\zc[t,t^{-1}]):A\equiv I \mod (t-\alpha)^i\}.$$
The sequence $\{K^i(\alpha)\}_{i\ge 1}$ is a central series in $K(\alpha) = K^1(\alpha)$
(i.e., $K^{i+j}(\alpha)\supseteq [K^i(\alpha),K^j(\alpha)]$).  Moreover, the graded
quotients satisfy
$$K^i(\alpha)/K^{i+1}(\alpha) \cong \lslnc.$$  The conjugation action of $\glnct$ on
$K(\alpha)$ induces a homomorphism
$$f_n(\alpha):\glnct\lra \text{Aut}(K(\alpha)/K^2(\alpha)) \cong GL_{n(n-2)}(\zc)$$
(note that $\lslnc$ is a vector space of dimension $(n-1)^2 - 1 = n(n-2)$).  We define
$$\rho_n(\alpha) = f_n(\alpha)\circ \beta_n.$$

The kernel of $\rho_n(1)$ is easily described:  it is the subgroup $P_n$ of pure braids.
The map $\rho_n(1)$ is thus a representation of the symmetric group $\Sigma_n$.

\medskip

\noindent {\bf Theorem \ref{rhothm}.}  {\em Suppose $n\ge 4$.  Denote by $V$ the standard
$(n-1)$-dimensional representation of $\Sigma_n$ and by $W$ the representation corresponding
to the partition $n-2,2$.  Then 
$$\rho_n(1) \cong V\oplus \bigwedge\nolimits^2 V\oplus W.$$
Also, $\rho_3(1)\cong V\oplus \bigwedge^2 V$.}

\medskip

When $\alpha\ne 1$, however, the image of $\rho_n(\alpha)$ is infinite.  Denote by $\galpha$ the subgroup
of $SL_{n(n-2)}(\zz[\alpha^{\pm 1}])$ consisting of matrices congruent to the identity modulo
$(\alpha - 1)$.  Then we have the following result.

\medskip

\noindent {\bf Proposition \ref{rhocongruence}.} {\em The image of $P_n$ under $\rho_n(\alpha)$ lies in
 $\galpha$.}

\medskip

For the representations $\mu_n(\alpha,\beta)$, we begin with the Lawrence--Krammer--Bigelow
representation
$$\kappa_n:B_n\lra \glnzq.$$  For each $i\ge 1$, define a subgroup $L^i(\alpha,\beta)$ by
$$L^i(\alpha,\beta) = \{A\in SL_{\binom{n}{2}}(\zc[t^{\pm 1},q^{\pm 1}]):A\equiv I \mod (t-\alpha,q-\beta)^i\}.$$
The first graded quotient satisfies
$$L(\alpha,\beta)/L^2(\alpha,\beta) \cong {\mathfrak sl}_{\binom{n}{2}}(\zc) \times {\mathfrak sl}_{\binom{n}{2}}(\zc).$$
The conjugation action of $\glncq$ on $L(\alpha,\beta)/L^2(\alpha,\beta)$ is diagonal; that is, if 
$(v,w)\in L(\alpha,\beta)/L^2(\alpha,\beta)$ then a matrix $A$ acts as
$$A:(v,w)\mapsto (AvA^{-1},AwA^{-1}).$$
Let $N=\binom{n}{2}^2 - 1$ and let $g_n(\alpha,\beta)$ be the homomorphism
$$g_n(\alpha,\beta):\glncq\lra GL_N(\zc)$$
obtained by restricting the conjugation action to the first factor of ${\mathfrak sl}_{\binom{n}{2}}(\zc)$.
Define $\mu_n(\alpha,\beta)$ to be the composite
$$\mu_n(\alpha,\beta)= g_n(\alpha,\beta)\circ \kappa_n.$$

The maps $\mu_n(\alpha,\beta)$ are more complicated than the $\rho_n(\alpha)$ mostly because the Burau
matrices $\beta_n(\sigma_i)$ of the braid generators are block diagonal, while the matrices
$\kappa_n(\sigma_i)$ are not.  We content ourselves to analyze a few special cases.  For example,
the kernel of $\mu_n(1,1)$ is $P_n$ and so $\mu_n(1,1)$ is a representation of $\Sigma_n$.  We give
the decomposition of $\mu_n(1,1)$ for $n=3,4,5$.

In Section \ref{faithful}, we discuss the faithfulness of the maps $\rho_n(\alpha)$ and $\mu_n(\alpha,\beta)$.
Of course, none of them is faithful since the center of $B_n$ lies in the kernel of each.  Moreover,
since $\beta_n$ is not faithful for $n\ge 5$, there must be additional elements in the kernel of $\rho_n(\alpha)$.
However, the map $\kappa_n$ is faithful for all $n$, and this allows us to deduce the following result.

\medskip

\noindent {\bf Theorem \ref{mufaith}.} {\em If $\alpha$ and $\beta$ are algebraically independent, then
the kernel of the representation $\mu_n(\alpha,\beta)$ is precisely the center of $B_n$.}

\medskip

In the particular case of $B_4$, it is known \cite{birman} that $\beta_4$ is faithful if and only if the
matrices $\beta_4(\sigma_3\sigma_1^{-1})$ and $\beta_4(\sigma_2\sigma_3\sigma_1^{-1}\sigma_2^{-1})$ generate a
free group of rank $2$.  We have the following result, which motivated the study of the $\rho_n(\alpha)$ in
the first place.

\medskip

\noindent {\bf Theorem \ref{powerfree}.} {\em There is a positive integer $M$ such that for any $m\ge M$, the
matrices $\beta_4(\sigma_3\sigma_1^{-1})^m$ and $\beta_4(\sigma_2\sigma_3\sigma_1^{-1}\sigma_2^{-1})^m$ generate
a free group of rank $2$.}

\medskip

This was proved first by S.~Moran \cite{moran}, who showed also that one can take $M=3$.  It was our hope
that passing to the map $\rho_4(\alpha)$ would allow us to show that $M=1$ is possible.  We show in Section
\ref{faithful} that the method of proof fails for $M=1,2$.

\medskip

\noindent {\em Acknowledgements.}  I thank Dan Cohen and Alex Suciu for inviting me to participate in the
Special Session on Arrangements in Topology and Algebraic Geometry at the AMS meeting in Baton Rouge, 
March, 2003.  The ensuing pressure to give a talk yielded these results.
  
\section{The representations $\rho_n(\alpha)$}\label{rho} 
Recall the reduced Burau representation $\beta_n:B_n\ra\glnzt$ defined as follows.  Let
$\sigma_1,\sigma_2,\dots ,\sigma_{n-1}$ be the standard generators of $B_n$ and set
$$\beta_n(\sigma_1)= \left[\begin{array}{rrcc}
                           -t & 1 & \cdots & 0 \\
                           0 & 1 & \cdots & 0 \\
                           \vdots & & \ddots & \\
                           0 & & \cdots & 1
                           \end{array}\right], \qquad \beta_n(\sigma_r)= \left[\begin{array}{rrcrccc}
                                                                      1 & & & & & & \\
                                                                        & \ddots & & & & & \\
                                                                        & & 1 & 0 & 0 & & \\
                                                                        & & t & -t & 1 & & \\
                                                                        & & 0 & 0 & 1 & & \\
                                                                        & & & & & \ddots & \\
                                                                        & & & & & & 1
                                                                        \end{array}\right]$$
where the row containing $t\;-t\;1$ in the matrix for $\sigma_r$ is the $r$th row, $1<r<n$.
The map $\beta_n$ is not faithful for $n\ge 5$ \cite{bigelow}.  Let $\alpha\in\zc^\times$ and
for each $i\ge 1$, define
$$K^i(\alpha)=\{A\in SL_{n-1}(\zc[t,t^{-1}]):A\equiv I\mod (t-\alpha)^i\}.$$  One checks easily
that $K^\bullet(\alpha)$ is a descending central series in $K(\alpha)=K^1(\alpha)$. If $A\in K^i(\alpha)$,
we may write $A = I + (t-\alpha)^i Y \mod (t-\alpha)^{i+1}$, where $Y$ is a matrix with entries in
$\zc$.  Since $\det A=1$, we have $\text{trace}(Y)=0$.  Define a map 
$$\pi:K^i(\alpha) \lra \lslnc$$ by $\pi(A) = Y$.  This is easily seen to be a surjective homomorphism
and the kernel of $\pi$ is the subgroup $K^{i+1}(\alpha)$.  Thus, $\pi$ induces an isomorphism
$$K^i(\alpha)/K^{i+1}(\alpha) \cong \lslnc.$$

Now let $i=1$.  Consider the following matrices in $K(\alpha)$:
\begin{eqnarray*}
A_{ij} & = & I + (t-\alpha)e_{ij}, \quad i\ne j,\, 1\le i,j\le n-1 \\
B_{ii} & = & I + (t-\alpha)e_{ii} - (t-\alpha)e_{i+1,i+1} - (t-\alpha)e_{i,i+1} + (t-\alpha) e_{i+1,i},\\
        &   & {} \qquad\qquad\qquad 1\le i\le n-2
\end{eqnarray*}
where $e_{ij}$ is the matrix having $1$ in the $i,j$ position and zeros elsewhere.  As a basis of
$K(\alpha)/K^2(\alpha)$, we choose the matrices $A_{ij}$, $1\le i,j\le n-1$ and $A_{ii} = B_{ii}A_{i,i+1}A_{i+1,i}^{-1}$,
$1\le i\le n-2$.  Note that under the map $\pi:K(\alpha)\ra\lslnc$, we have
$$\pi(A_{ij}) = e_{ij},\quad\text{and}\quad \pi(A_{ii}) = e_{ii} - e_{i+1,i+1}.$$

The group $\glnct$ acts on $K(\alpha)$ via conjugation and hence acts on the quotient $K(\alpha)/K^2(\alpha)$.
Denote by $f_n(\alpha)$ the map
$$f_n(\alpha):\glnct\lra \text{Aut}(K(\alpha)/K^2(\alpha)) \cong GL_{n(n-2)}(\zc)$$
(note that $\dim \lslnc = (n-1)^2 - 1 = n(n-2)$).  Restricting this action to the $\beta_n(\sigma_i)$ gives a map
$$\rho_n(\alpha):B_n\lra GL_{n(n-2)}(\zc)$$
(that is, $\rho_n(\alpha) = f_n(\alpha)\circ\beta_n$).  The action of each $\sigma_i$ on $K(\alpha)/K^2(\alpha)$
is given by the following formul\ae:
$$\begin{array}{rll}
\sigma_1: & A_{ij}\mapsto A_{ij} & 3\le i\le n-1 \\
          &                      & 2\le j\le n-1 \\
          & A_{ii}\mapsto A_{ii} & 3\le i\le n-2 \\
          & A_{1j}\mapsto -\alpha A_{1j} & 3\le j\le n-1 \\
          & A_{2j}\mapsto A_{1j} + A_{2j} & 3\le j\le n-1 \\
          & A_{i1}\mapsto -\frac{1}{\alpha}A_{i1} + \frac{1}{\alpha}A_{i2} & 3\le i\le n-1 \\
          & A_{12}\mapsto -\alpha A_{12} & \\
          & A_{21}\mapsto \frac{1}{\alpha}A_{12} - \frac{1}{\alpha}A_{21}-\frac{1}{\alpha}A_{11} & \\
          & A_{11}\mapsto A_{11}-2A_{12} &            \\
          & A_{22}\mapsto A_{12}+A_{22}       & \\
\end{array}$$
and for $2\le k\le n-2$ and $i,j\ne k-1,k,k+1$
$$\begin{array}{rl}
\sigma_k: & A_{ij}\mapsto A_{ij}  \\
        & A_{ii}\mapsto A_{ii} \\
          & A_{i,k-1}\mapsto A_{i,k-1} \\
          & A_{ik}\mapsto A_{i,k-1} -\frac{1}{\alpha}A_{ik} +\frac{1}{\alpha}A_{i,k+1}  \\
          & A_{i,k+1}\mapsto A_{i,k+1} \\
          & A_{k-1,j}\mapsto A_{k-1,j}+\alpha A_{kj} \\
          & A_{kj}\mapsto -\alpha A_{kj}  \\
          & A_{k+1,j}\mapsto A_{kj}+A_{k+1,j}  \\
          & A_{k-1,k} \mapsto A_{k-1,k-1}+ -\frac{1}{\alpha}A_{k-1,k}+\frac{1}{\alpha}A_{k-1,k+1}+\alpha A_{k,k-1}+A_{k,k+1} \\
          & A_{k-1,k+1}\mapsto A_{k-1,k+1}+\alpha A_{k,k+1} \\
          & A_{k,k-1}\mapsto -\alpha A_{k,k-1}  \\
          & A_{k,k+1}\mapsto -\alpha A_{k,k+1}  \\
          & A_{k+1,k-1}\mapsto A_{k,k-1}+A_{k+1,k-1}  \\
          & A_{k+1,k}\mapsto -\frac{1}{\alpha}A_{kk} +A_{k,k-1}+A_{k+1,k-1}+\frac{1}{\alpha}A_{k,k+1}-\frac{1}{\alpha}A_{k+1,k}  \\
          & A_{k-2,k-2}\mapsto A_{k-2,k-2}-\alpha A_{k,k-1} \quad (k\ge 3) \\
          & A_{k-1,k-1}\mapsto A_{k-1,k-1} + 2\alpha A_{k,k-1} + A_{k,k+1}  \\
          & A_{kk}\mapsto A_{kk} - \alpha A_{k,k-1} - 2A_{k,k+1} \\
          & A_{k+1,k+1}\mapsto A_{k+1,k+1} +A_{k,k+1} \quad (k\le n-3) \\
\end{array}$$
and finally
$$\begin{array}{rll}
\sigma_{n-1}: & A_{ij}\mapsto A_{ij} & 1\le i\le n-3 \\
               &                     &  1\le j\le n-2 \\
              & A_{ii}\mapsto A_{ii} & 1\le i\le n-4 \\
              & A_{n-2,j}\mapsto A_{n-2,j}+\alpha A_{n-1,j} & 1\le j\le n-3 \\
              & A_{n-1,j}\mapsto -\alpha A_{n-1,j} & 1\le j\le n-2 \\
              & A_{i,n-1}\mapsto A_{i,n-2}-\frac{1}{\alpha}A_{i,n-1} & 1\le i\le n-3 \\
              & A_{n-2,n-1}\mapsto -\frac{1}{\alpha}A_{n-2,n-1} + \alpha A_{n-1,n-2} +A_{n-2,n-2} & \\
              & A_{n-3,n-3}\mapsto  A_{n-3,n-3}-\alpha A_{n-1,n-2} & \\
              & A_{n-2,n-2}\mapsto 2\alpha A_{n-1,n-2} + A_{n-2,n-2} & \\
\end{array}$$

Consider the particular case $n=3$.  We have $\dim \rho_3(\alpha)=3$.  The $3$-dimensional representations
of $B_3$ were classified by Tuba and Wenzl \cite{tuba}---they are characterized uniquely be the eigenvalues
of the matrices for $\sigma_1$ and $\sigma_2$.  Using the basis $e_1=-A_{12}$, $e_2=A_{11}$,
$e_3 = -\alpha A_{21}$, the matrices of $\rho_3(\alpha)(\sigma_1)$ and $\rho_3(\alpha)(\sigma_2)$ are
$$\rho_3(\alpha)(\sigma_1) = \left[\begin{array}{rrr}
                                 -\alpha & 2 & 1 \\
                                    0 & 1 & 1 \\
                                    0 & 0 & -\frac{1}{\alpha}
                                    \end{array}\right],\qquad \rho_3(\alpha)(\sigma_2) = \left[\begin{array}{rrr}
                                                                                             -\frac{1}{\alpha} & 0 & 0 \\
                                                                                             -1 & 1 & 0 \\
                                                                                             1 & -2 & -\alpha
                                                                                             \end{array}\right].$$
These correspond to the matrices in \cite{tuba} with $\lambda_1=-\alpha$, $\lambda_2=1$ and $\lambda_3=-\frac{1}{\alpha}$.
We have thus found a naturally occurring one-parameter family of the abstract representations defined in \cite{tuba}.

Observe that the map $\rho_n(\alpha)$ is defined over the ring $\zz[\alpha,\alpha^{-1}]$.  Denote by $\galpha$ the subgroup
of $GL_{n(n-2)}(\zz[\alpha,\alpha^{-1}])$ consisting of matrices congruent to $I$ modulo $(\alpha - 1)$.  Let $P_n$ be
the subgroup of pure braids.

\begin{prop}\label{rhocongruence}  The image of $P_n$ under $\rho_n(\alpha)$ lies in $\galpha$.
\end{prop}

\begin{proof} The group $P_n$ is generated by the $\sigma_i^2$ along with some conjugates of these.  Since
the $\sigma_i$ are all conjugate in $B_n$, so are the $\sigma_i^2$.  Thus it suffices to show that
$\rho_n(\alpha)(\sigma_1^2)\in\galpha$.  But this is easy:
$$\begin{array}{rll}
\sigma_1^2: & A_{ij}\mapsto A_{ij} & 3\le i,j\le n-1 \\
            & A_{ii}\mapsto A_{ii} & 3\le i\le n-2 \\
            & A_{1,j}\mapsto (1+2(\alpha-1)+(\alpha-1)^2)A_{1,j} & 2\le j\le n-1 \\
            & A_{2,j}\mapsto (1-\alpha)A_{1,j} + A_{2,j} & 3\le j\le n-1 \\
            & A_{i,1}\mapsto (1 + \frac{1-\alpha}{\alpha}+ \frac{1-\alpha}{\alpha^2})A_{i,1} + \frac{\alpha-1}{\alpha^2}A_{i,2} & 3\le i\le n-1 \\
            & A_{i,2}\mapsto A_{i,2} & 3\le i\le n-1 \\
            & A_{21}\mapsto -\frac{(1-\alpha)^2}{\alpha}A_{12} + (1 + \frac{1-\alpha}{\alpha}+ \frac{1-\alpha}{\alpha^2})A_{21} + \frac{1-\alpha}{\alpha^2}A_{11} & \\
            & A_{11}\mapsto A_{11}-2(1-\alpha)A_{12} & \\
            & A_{22}\mapsto (1-\alpha)A_{12} + A_{22} & 
\end{array}$$
\end{proof}

In particular, if $\alpha=1$, then $P_n\subseteq \text{ker} \rho_n(1)$.  (This also follows from the fact that
$\beta_n(P_n)\subset K(1)$ and hence the $\beta_n(\sigma_i^2)$ act trivially on $K(1)/K^2(1)$.)  The reverse inclusion
also holds since the action of any $\rho_n(1)(\sigma_i)$ is the same as the action of $\beta_n(\sigma_i)$ evaluated
at $t=1$, and this collection of matrices is a faithful representation of the symmetric group $\Sigma_n$.  For example,
when $n=3$ we obtain (using the basis $e_1=A_{12}$, $e_2=A_{21}$, $e_3=A_{11}-A_{12}+A_{21}$)
$$\rho_3(1)(\sigma_1) = \left[\begin{array}{rrr}
                             -1 & 0 & 0 \\
                             0 & 0 & -1 \\
                             0 & -1 & 0
                             \end{array}\right]\qquad \rho_3(1)(\sigma_2) = \left[\begin{array}{rrr}
                                                                                0 & 0 & -1 \\
                                                                                0 & -1 & 0 \\
                                                                                -1 & 0 & 0
                                                                                \end{array}\right].$$
The character of this is
$$\begin{array}{rccc}
      &  (1) & (12) & (123) \\ \hline
\chi_{\rho_3(1)} & 3 & -1 & 0
\end{array}$$
If $V$ denotes the standard $2$-dimensional representation of $\Sigma_3$, then
$\chi_{\rho_3(1)}=\chi_V + \chi_{\bigwedge^2V}$.  Thus, $\rho_3(1) \cong V\oplus \bigwedge^2 V$.

For $n=4$ we have
$$\begin{array}{rccccc}
            & (1) & (12) & (123) & (1234) & (12)(34) \\ \hline
\chi_{\rho_4(1)} & 8 & 0 & -1 & 0 & 0
\end{array}$$
Denote by $W$ the representation corresponding to the partition $2,2$.  Then an easy check shows
that $\chi_{\rho_4} = \chi_V + \chi_{\bigwedge^2 V} + \chi_W$.

The character of $\rho_5(1)$ is
$$\begin{array}{cccccccc}
    & (1) & (12) & (123) & (1234) & (12345) & (12)(34) & (12)(345) \\ \hline
\chi_{\rho_5(1)} & 15 & 3 & 0 & -1 & 0 & -1 & 0
\end{array}$$

\begin{theorem}\label{rhothm}  Suppose $n\ge 4$.  Denote by $V$ the standard $(n-1)$-dimensional
representation of $\Sigma_n$ and by $W$ the representation corresponding to the partition $n-2,2$. Then
$$\rho_n(1) \cong V\oplus \bigwedge\nolimits^2 V\oplus W.$$
\end{theorem}

\begin{proof}  First note that $\dim V = n-1$, $\dim \bigwedge^2 V = \frac{(n-1)(n-2)}{2}$ and $\dim W
=\frac{n(n-3)}{2}$ (\cite{fulton}, p.~50).  Thus,
$$(n-1) + \frac{(n-1)(n-2)}{2} + \frac{n(n-3)}{2} = n(n-2)$$
so that $\dim \rho_n(1) = \dim(V\oplus\bigwedge^2 V\oplus W)$.
Now if $C_{\underline{i}}$, $\underline{i} = (i_1,i_2,\dots ,i_n)$ denotes the conjugacy class in $\Sigma_n$ 
consisting of cycles with $i_1$ $1$-cycles, $i_2$ 2-cycles, etc., then
$$(\chi_V + \chi_{\bigwedge\nolimits^2 V} + \chi_W)(C_{\underline{i}}) = i_1(i_1-2)$$
(this follows from \cite{fulton}, 4.15, p.~51).  Direct calculation shows that
\begin{eqnarray*}
\chi_{\rho_n(1)}((12)) & = & (n-2)(n-4) = i_1(i_1-2)   \\
\chi_{\rho_n(1)}((123)) & = & (n-3)(n-5) = i_1(i_1-2)  \\
\chi_{\rho_n(1)}((12)(n-1,n)) & = & (n-4)(n-6) = i_1(i_1-2),
\end{eqnarray*}
etc.  Thus, $\chi_{\rho_n(1)} = \chi_V + \chi_{\bigwedge^2 V} + \chi_W$ and since $V$, $\bigwedge^2 V$, and
$W$ are distinct irreducible representations, the result follows.
\end{proof}

\section{The representations $\mu_n(\alpha,\beta)$}\label{mu}
Let $R=\zz[t^{\pm 1},q^{\pm 1}]$ be the ring of Laurent polynomials in $t,q$ and let $A$ be the
free $R$-module
$$A=\bigoplus_{1\le i<j\le n} Rx_{ij}.$$  The Lawrence--Krammer--Bigelow representation
$$\kappa_n:B_n\lra \glnzq$$ is defined by
$$\begin{array}{rll}
\kappa_n(\sigma_k): & x_{k,k+1}\mapsto tq^2x_{k,k+1} & \\
            & x_{ik}\mapsto (1-q)x_{ik} + qx_{i,k+1} & i<k \\
            & x_{i,k+1}\mapsto x_{ik}+tq^{k-i+1}(q-1)x_{k,k+1} & i<k \\
            & x_{kj}\mapsto tq(q-1)x_{k,k+1}+qx_{k+1,j} & k+1<j \\
            & x_{k+1,j}\mapsto x_{kj} + (1-q)x_{k+1,j} & k+1<j \\
            & x_{ij}\mapsto x_{ij} & i<j<k\;\text{or}\; k+1<i<j \\
            & x_{ij}\mapsto x_{ij}+tq^{k-i}(q-1)^2x_{k,k+1} & i<k<k+1<j.
\end{array}$$
The map $\kappa_n$ is faithful for all $n$ \cite{bigelow2},\cite{krammer}.  This is the only known
faithful representation of $B_n$, $n\ge 4$.

For each $i\ge 1$ define, for $\alpha,\beta\in\zc^\times$,
$$L^i(\alpha,\beta) = \{A\in SL_{\binom{n}{2}}(\zc[t^{\pm 1},q^{\pm 1}]): A\equiv I\mod (t-\alpha,q-\beta)^i\}$$
(here, $(t-\alpha,q-\beta)$ denotes the ideal generated by $t-\alpha$ and $q-\beta$).
This is a descending central series in $L(\alpha,\beta)=L^1(\alpha,\beta)$.  The graded quotients are
more complicated here, but we are interested only in the first one: $L(\alpha,\beta)/L^2(\alpha,\beta)$.  If
$A\in L(\alpha,\beta)$, write
$$A=I+(t-\alpha)A_t + (q-\beta)A_q + X,$$  where $A_t,A_q$ are $\binom{n}{2}\times\binom{n}{2}$ matrices
over $\zc$ and $X\equiv 0 \mod (t-\alpha,q-\beta)^2$.  Again, the condition $\det A = 1$ forces
$\text{tr}(A_t)=0=\text{tr}(A_q)$.  Define a map
$$\pi:L(\alpha,\beta)\lra {\mathfrak sl}_{\binom{n}{2}}(\zc)\times {\mathfrak sl}_{\binom{n}{2}}(\zc)$$
by $\pi(A) = (A_t,A_q)$.  This is a surjective group homomorphism with kernel $L^2(\alpha,\beta)$.

If $Z\in\glncq$, we may write $Z=Z_0Z_1$, where $Z_0\in GL_{\binom{n}{2}}(\zc)$ and $Z_1\in L(\alpha,\beta)$.
Then $Z_1$ acts trivially on $L(\alpha,\beta)/L^2(\alpha,\beta)$ and so $Z$ acts on $L(\alpha,\beta)/L^2(\alpha,\beta)$
by
$$Z:(A_t,A_q)\mapsto (Z_0A_tZ_0^{-1},Z_0A_qZ_0^{-1});$$
that is, the action is diagonal.  Consider the action on the first factor (call it $L$):
$$g_n(\alpha,\beta):\glncq\lra \text{Aut}(L)\cong GL_N(\zc)$$
where $N=\binom{n}{2}^2-1$.  Restricting this to the image of $B_n$ under $\kappa_n$ yields a map
$$\mu_n(\alpha,\beta):B_n\lra GL_N(\zc)$$
(i.e., $\mu_n(\alpha,\beta) = g_n(\alpha,\beta)\circ\kappa_n$).

We shall not write down a formula for the $\mu_n(\alpha,\beta)$ in general.  Indeed, the following
formula for $\mu_3(\alpha,\beta)$ shows that the general case is hopelessly complicated. Using the basis
described in (\ref{basis1}) and (\ref{basis2}) below, the matrix
of $\mu_3(\alpha,\beta)(\sigma_1)$ is
{\Tiny
$$\left[\begin{array}{cccccccc}
\alpha\beta(\beta-1) & \alpha\beta^2 & -\frac{\alpha}{\beta}(\beta-1)^3 & \alpha\beta(\beta-1) & 0 & 0 & -2\alpha(\beta-1)^2 & \alpha(\beta-1)^2 \\
\alpha\beta & 0 & -\frac{\alpha}{\beta}(\beta-1)^2 & 0 & 0 & 0 & -2\alpha(\beta-1) & \alpha(\beta-1) \\
0 & 0 & 0 & 0 & \frac{\alpha}{\beta^{2}} & 0 & 0 & 0 \\
0 & 0 & 0 & 0 & -\frac{(\beta-1)}{\beta} & \frac{1}{\beta} & 0 & 0 \\
0 & 0 & \frac{1}{\alpha\beta} & 0 & -\frac{\beta-1}{\alpha\beta^{2}} & 0 & 0 & 0 \\
0 & 0 & -\frac{(\beta-1)^2}{\beta} & \beta & \frac{(\beta-1)^3}{\beta^2} & -\frac{(\beta-1)^2}{\beta} & 1-\beta & 2(\beta-1) \\
0 & 0 & 1-\frac{1}{\beta} & 0 & 0 & 0 & 1 & 0 \\
0 & 0 & 1-\frac{1}{\beta} & 0 & -\frac{(\beta-1)^2}{\beta^2} & 1-\frac{1}{\beta} & 1 & -1
\end{array}\right]
$$}
and that of $\mu_3(\alpha,\beta)(\sigma_2)$ is
{\Tiny
$$\left[\begin{array}{cccccccc}
-\frac{(\beta-1)^2}{\beta} & \frac{(\beta-1)^3}{\beta} & \frac{1}{\beta} & -\frac{(\beta-1)^2}{\beta} & 0 & 0 & \frac{-2(\beta-1)}{\beta} & \frac{\beta-1}{\beta} \\
0 & -\frac{\beta-1}{\alpha\beta^2} & 0 & \frac{1}{\alpha\beta^2} & 0 & 0 & 0 & 0 \\
\beta & -\beta(\beta-1) & 0 & 0 & 0 & 0 & 0 & 0 \\
0 & \frac{1}{\alpha\beta} & 0 & 0 & 0 & 0 & 0 & 0 \\
0 & 0 & 0 & -\alpha\beta^2(\beta-1)^2 & 0 & \alpha\beta^2 & -\alpha\beta^2(\beta-1) & 2\alpha\beta^2(\beta-1) \\
0 & 0 & \alpha\beta(\beta-1) & -\alpha\beta(\beta-1)^3 & \alpha\beta & \alpha\beta(\beta-1) & -\alpha\beta(\beta-1)^2 & 2\alpha\beta(\beta-1)^2 \\
1-\beta & (\beta-1)^2 & 0 & 1-\beta & 0 & 0 & -1 & 1 \\
0 & 0 & 0 & 1-\beta & 0 & 0 & 0 & 1
\end{array}\right].
$$}
Note however that $\mu_n(\alpha,\beta)$ is defined over $\zz[\alpha^{\pm 1},\beta^{\pm 1}]$ and we have the following
result.  Denote by $\Gamma_N(\alpha,\beta)$ the subgroup of $GL_N(\zz[\alpha^{\pm 1},\beta^{\pm 1}])$ consisting of
those matrices that are congruent to $I$ modulo $(\alpha-1,\beta-1)$.

\begin{prop}\label{mucongruence} The image of $P_n$ under $\mu_n(\alpha,\beta)$ lies in $\Gamma_N(\alpha,\beta)$.
\end{prop}

\begin{proof} It suffices to check this for the $\mu_n(\alpha,\beta)(\sigma_k^2)$.  First note the following
formula for $\kappa_n(\sigma_k^2)$:
$$\begin{array}{l}
  x_{k,k+1}\mapsto t^2q^4x_{k,k+1}  \\
  x_{ik}\mapsto (1+(q-1)+(q-1)^2)x_{ik}+(1-q)qx_{i,k+1}+tq^{k-i+1}(q-1)x_{k,k+1}  \\
  x_{i,k+1}\mapsto (1-q)x_{ik}+(1+(q-1))x_{i,k+1}+t^2q^{k-i+3}(q-1)x_{k,k+1}  \\
  x_{kj}\mapsto t^2q^3(q-1)x_{k,k+1} + (1+(q-1))x_{kj}+q(1-q)x_{k+1,j}  \\
  x_{k+1,j}\mapsto tq(q-1)x_{k,k+1} + (1+(q-1)+(q-1)^2)x_{k+1,j}+(1-q)x_{kj} \\
  x_{ij}\mapsto x_{ij} + tq^{k-i}(q-1)^2(1+tq^2)x_{k,k+1} 
\end{array}$$
where the ranges on $i,j$ are ($i<k$), ($i<k$), ($k+1<j$), ($k+1<j$), ($i<j<k\;\text{or}\;k+1<i<j$),
and ($i<k<k+1<j$), respectively.  Note that $t^2q^4\equiv 1+2(t-1) + 4(q-1) \mod (t-1,q-1)^2$.  Thus,
after evaluating $\kappa_n(\sigma_k^2)$ at $t=\alpha$, $q=\beta$, we may write
$$\kappa_n(\sigma_k^2) = X_kY_k,$$
where $X_k\in L(\alpha,\beta)$ and $Y_k\equiv I \mod (\alpha-1,\beta-1)$.  The action of $\kappa_n(\sigma_k^2)$
on $L$ is via conjugation by $Y_k$.  Write
$$Y_k=I + (\alpha-1)Y_\alpha + (\beta-1)Y_\beta + Z$$
where $Z\equiv 0\mod (\alpha-1,\beta-1)^2$.  Then if $A=I+(t-\alpha)A_t + X$ ($X\equiv 0\mod (t-\alpha,q-\beta)^2$),
we have
{\tiny
\begin{eqnarray*}
Y_kAY_k^{-1} & = & (1+(\alpha-1)Y_\alpha+(\beta-1)Y_\beta+Z)(I+(t-\alpha)A_t+X)(I-(\alpha-1)Y_\alpha-(\beta-1)Y_\beta+U) \\
             & \equiv & I +(t-\alpha)A_t + (\alpha-1)(t-\alpha)[Y_\alpha,A_t] + (\beta-1)(t-\alpha)[Y_\beta,A_t] \mod (t-\alpha,q-\beta)^2 \\
             & \equiv & A \mod (\alpha-1,\beta-1).
\end{eqnarray*}}
Thus, $\kappa_n(\sigma_k^2)$ acts as the identity on $L$ modulo $(\alpha-1,\beta-1)$; that is,
$\mu_n(\alpha,\beta)(\sigma_k^2)\in \Gamma_N(\alpha,\beta)$.
\end{proof}

\begin{cor}\label{musigmarep} The kernel of $\mu_n(1,1)$ is the subgroup $P_n$ and so $\mu_n(1,1)$
is a representation of $\Sigma_n$.\hfill $\qed$
\end{cor}

As a basis of $L$ we take the matrices
\begin{equation}\label{basis1}
A_{ij}=I+(t-\alpha)e_{ij} \qquad 1\le i,j \le \binom{n}{2}, i\ne j,
\end{equation}
and
\begin{equation}\label{basis2}
A_{ii} \equiv I+(t-\alpha)e_{ii}-(t-\alpha)e_{i+1,i+1} \qquad 1\le i\le \binom{n}{2}-1,
\end{equation}
ordered as $A_{12},A_{13},\dots ,A_{1,n},A_{21},\dots ,A_{n,n-1},A_{11},A_{22},\dots$.
It is possible to give a simple formula for $\mu_n(1,1)$.  Note that upon evaluating
$\kappa_n(\sigma_k)$ at $t=1,q=1$, one obtains the permutation matrix
$$\begin{array}{rll}
\tau_k: & x_{k,k+1}\mapsto x_{k,k+1} & \\
        & x_{ik}\mapsto x_{i,k+1} & i<k \\
        & x_{i,k+1}\mapsto x_{ik} & i<k \\
        & x_{kj}\mapsto x_{k+1,j} & k+1<j \\
        & x_{k+1,j}\mapsto x_{kj} & k+1<j \\
        & x_{ij}\mapsto x_{ij} & i<j<k,k+1<i<j,i<k<k+1<j.
\end{array}$$
Order the basis of $A=\displaystyle\bigoplus_{1\le i<j\le n}Rx_{ij}$ as
$$e_1=x_{12}, e_2=x_{13},\dots ,e_{n-1}=x_{1,n},e_n=x_{23},\dots ,e_{\binom{n}{2}}=x_{n-1,n}$$
and consider $\tau_k$ as a permutation of the set $\{1,2,\dots ,\binom{n}{2}\}$.  The action of
$\mu_n(1,1)$ may then be described as
$$\mu_n(1,1)(\sigma_k): \begin{cases}
                           A_{ij}\mapsto A_{\tau_k(i),\tau_k(j)} & i\ne j \\
                                {}                                 &        \\
                           {A_{ii}\mapsto \begin{cases}
                                         \displaystyle\sum_{\ell=\tau_k(i)}^{\tau_k(i+1)-1} A_{\ell\ell} & \tau_k(i)<\tau_k(i+1) \\
                                         \displaystyle -\sum_{\ell=\tau_k(i+1)}^{\tau_k(i)-1} A_{\ell\ell} & \tau_k(i+1)<\tau_k(i).
                                         \end{cases}}
                           \end{cases}$$ 

The characters of $\mu_n(1,1)$ for $n=3,4,5$ are as follows:
$$\begin{array}{r|ccc}
                 & (1) & (12) & (123) \\ \hline
\chi_{\mu_3(1,1)} & 8 & 0 & -1
\end{array}$$
$$\begin{array}{r|ccccc}
                  & (1) & (12) & (123) & (1234) & (12)(34) \\ \hline
\chi_{\mu_4(1,1)} & 35 & 3 & -1 & 3 & -1 
\end{array}$$
$$\begin{array}{r|ccccccc}
                   & (1) & (12) & (123) & (1234) & (12345) & (12)(34) & (12)(345) \\ \hline
\chi_{\mu_5(1,1)} & 99 & 15 & 0 & -1 & -1 & 3 & 0
\end{array}$$
A straightforward calculation then shows the following.
{\small
\begin{eqnarray*}
\mu_3(1,1) & \cong & (\text{alt})\oplus (\text{triv})\oplus V^{\oplus 3} \\
\mu_4(1,1) & \cong & (\text{alt})\oplus (\text{triv})^{\oplus 2} \oplus (\bigwedge\nolimits^2V)^{\oplus 3}\oplus W^{\oplus 4}\oplus V^{\oplus 5} \\
\mu_5(1,1) & \cong & (V\otimes(\text{alt}))\oplus (\text{triv})^{\oplus 2}\oplus (W\otimes(\text{alt}))^{\oplus 3}\oplus (\bigwedge\nolimits^2V)^{\oplus 4}
\oplus W^{\oplus 6} \oplus V^{\oplus 6}
\end{eqnarray*}}
(recall that $V$ is the standard $(n-1)$-dimensional representation and $W$ is the representation corresponding to the
partition $n-2,2$).

Note that the $A_{ii}$, $1\le i\le \binom{n}{2}-1$ form a $\Sigma_n$-submodule of $L$.

\begin{prop}\label{submodule} The submodule $U$ spanned by the $A_{ii}$, $1\le i\le \binom{n}{2}-1$, is isomorphic
to $V\oplus W$.
\end{prop}

\begin{proof}  Recall that $W$ corresponds to the partition $n-2,2$; we have $\dim W=\displaystyle\frac{n(n-3)}{2}$.  Then
\begin{eqnarray*}
\dim V + \dim W & = & (n-1) + \frac{n(n-3)}{2} \\
                & = & \frac{n^2-n-2}{2} \\
                & = & \binom{n}{2} - 1 \\
                & = & \dim U.
\end{eqnarray*}
Note that $\tau_1=(2,n)(3,n+1)\cdots (n-1,2n-3)$ and so
\begin{eqnarray*}
\chi_U(12) = \text{tr}(\tau_1) & = & \sum_{j=2n-3}^{\binom{n}{2}-1} 1 \\
                               & = & \biggl(\binom{n}{2}-1\biggr)-(2n-3)+1 \\
                               & = & (n-3) + \frac{1}{2}(n-3)(n-4) \\
                               & = & (\chi_V+\chi_W)(12).
\end{eqnarray*}
The values of $\chi_U$ for other conjugacy classes are obtained similarly.
\end{proof}
 
\section{Faithfulness}\label{faithful}
It is clear that none of the representations $\rho_n(\alpha)$, $\mu_n(\alpha,\beta)$ is faithful.
Indeed, the center of $B_n$ lies in the kernel of each $\rho_n(\alpha)$ and $\mu_n(\alpha,\beta)$.

Is there more in the kernel?  The answer is certainly yes for $\rho_n(\alpha)$, $n\ge 5$ since
$\beta_n$ is not faithful.

\subsection{Associated graded algebras}\label{graded}  One approach is to study the map on associated
graded algebras.  Let $\alpha=-1$.  Then the image of $P_n$ under $\rho_n(-1)$ lies in the subgroup
$$\Gamma_{n(n-2)}(-1) = \{A\in SL_{n(n-2)}(\zz):A\equiv I\mod 2\}.$$
The lower central series of $\Gamma_{n(n-2)}(-1)$ is well-understood via the work of Bass--Milnor--Serre
\cite{bms}; the $i$th term of the lower central series is
$$\Gamma^i_{n(n-2)}(-1) = \{A\in\Gamma_{n(n-2)}(-1): A\equiv I\mod 2^i\}$$
and the graded quotients satisfy
$$\Gamma^i/\Gamma^{i+1}\cong {\mathfrak sl}_{n(n-2)}({\mathbb F}_2).$$

The structure of $\text{Gr}^\bullet P_n$ is known thanks to the work of Kohno \cite{kohno}.  Each
graded quotient $\Gamma^i P_n/\Gamma^{i+1} P_n$ is free abelian with rank $\varphi_i(n)$ given by
the formula
$$\prod_{i=1}^\infty (1-t^i)^{\varphi_i(n)} = \prod_{j=1}^{n-1} (1-jt).$$

Consider the map of associated graded algebras
$$\text{Gr}^\bullet \rho_n(-1):\text{Gr}^\bullet P_n\lra \text{Gr}^\bullet \Gamma_{n(n-2)}(-1).$$
Let us examine first the case $n=3$.  Here, using the basis $e_1=-A_{12}$, $e_2=A_{11}$ and
$e_3=A_{21}$, we have
$$\rho_3(-1)(\sigma_1) = \left[\begin{array}{ccc}
                                  1 & 2 & 1 \\
                                  0 & 1 & 1 \\
                                  0 & 0 & 1
                                  \end{array}\right]\qquad \rho_3(-1)(\sigma_2) = \left[\begin{array}{rrr}
                                                                                       1 & 0 & 0 \\
                                                                                       -1 & 1 & 0 \\
                                                                                       1 & -2 & 1
                                                                                       \end{array}\right].$$
Denote the generators of $P_3$ by
$$B_{12}=\sigma_1^2, \quad B_{13}=\sigma_2\sigma_1^2\sigma_2^{-1}, \quad B_{23} = \sigma_2^2.$$
Then the map
$$\text{Gr}^1\rho_3(-1): H_1(P_3;\zz)\lra H_1(\Gamma_3(-1);\zz)$$
is the map $\zz\{B_{12},B_{13},B_{23}\}\ra {\mathfrak sl}_3({\mathbb F}_2)$
\begin{eqnarray*}
B_{12} & \mapsto & e_{23} \\
B_{13} & \mapsto & e_{21}+e_{23} \\
B_{23} & \mapsto & e_{21}.
\end{eqnarray*}
Note that in ${\mathfrak sl}_3({\mathbb F}_2)$, $[e_{21},e_{23}]=0$ and so the image of $\text{Gr}^\bullet \rho_3(-1)$
is simply the submodule of $\text{Gr}^\bullet \Gamma_3(-1)$ spanned by $e_{12},e_{23}\in\text{Gr}^1\Gamma_3(-1)$; that is
$$\text{Gr}^i\rho_3(-1):\Gamma^i P_n/\Gamma^{i+1} P_n \lra \Gamma^i_3(-1)/\Gamma^{i+1}_3(-1)$$ is the zero map for $i\ge 2$.  In particular,
this tells us that if $x\in \Gamma^i P_3$, $i\ge 2$, then $\rho_3(-1)(x)$ is congruent to the identity matrix modulo
$2^{i+1}$ instead of $2^i$.

By contrast, for $n\ge 4$ the map $\text{Gr}^\bullet P_n\ra \text{Gr}^\bullet \Gamma_{n(n-2)}(-1)$ is highly
nontrivial.  Moreover, for $i$ large, the map
$$(\Gamma^i P_n/\Gamma^{i+1} P_n)\otimes {\mathbb F}_2\lra \Gamma^i_{n(n-2)}(-1)/\Gamma^{i+1}_{n(n-2)}(-1)$$
cannot be injective (the rank of the domain is greater than $(n(n-2))^2-1$ for $i$ large).  This gives a method
for searching for elements in the kernel of $\beta_n$---find an element in the kernel of $\text{Gr}^i\rho_n(-1)$,
lift it to $P_n$, and compute its Burau matrix.  Of course, this is terribly inefficient.

\subsection{The case of $B_4$}\label{burau4}
Let us examine the maps $\rho_4(\alpha)$ in greater detail.  According to Theorem 3.19 of \cite{birman},
$\beta_4$ is faithful if and only if the matrices
$$x=\beta_4(\sigma_3\sigma_1^{-1})\qquad \text{and}\qquad y=\beta_4(\sigma_2\sigma_3\sigma_1^{-1}\sigma_2^{-1})$$
generate a rank 2 free subgroup of $GL_3(\zz[t,t^{-1}])$.  In turn, this will hold if and only if
the matrices
$$X(\alpha)=\rho_4(\alpha)(\sigma_3\sigma_1^{-1})\qquad \text{and}\qquad Y(\alpha)=\rho_4(\alpha)(\sigma_2\sigma_3\sigma_1^{-1}\sigma_2^{-1})$$
($\alpha\ne 1$) generate a rank 2 free subgroup of $SL_8(\zc)$.  We do not have a proof of this, but we
do have the following.

\begin{theorem}\label{powerfree} There is a positive integer $M$ such that for all $m\ge M$, the group generated
by $X(\alpha)^m$ and $Y(\alpha)^m$ is free of rank $2$.
\end{theorem}

\begin{proof}  We use Proposition 3.12 of \cite{tits}.  We first establish notation.  If $g\in GL(V)$, where
$V$ is a complex vector space, write its characteristic polynomial as $\prod_{i=1}^n(t-\alpha_i)$ and let
$\Omega = \{\alpha_i:\alpha_i=\max_{j}\{||\alpha_j||\}\}$.  Define polynomials $f_1(t)$ and $f_2(t)$ by
$$f_1(t)=\prod_{\alpha\in\Omega}(t-\alpha)\qquad \text{and}\qquad f_2(t)=\prod_{\alpha\not\in\Omega}(t-\alpha).$$
Let $A(g)$ be the subspace of ${\mathbb P}_V$ (the projective space associated to $V$) corresponding
to the kernel of $f_1(g)$ and let $A'(g)$ be that corresponding to the kernel of $f_2(g)$.

To prove the theorem, we need only show
\begin{eqnarray}
 {} & & A(X(\alpha)), A(X(\alpha)^{-1}), A(Y(\alpha)), A(Y(\alpha)^{-1})\; \text{are points}  \\ \label{cond2}
 {} & & A(X(\alpha))\cup A(X(\alpha)^{-1}) \subset {\mathbb P}_V - A'(Y(\alpha))-A'(Y(\alpha)^{-1}) \\ \label{cond3} 
 {} & & A(Y(\alpha))\cup A(Y(\alpha)^{-1}) \subset {\mathbb P}_V - A'(X(\alpha))-A'(X(\alpha)^{-1}). 
\end{eqnarray}
The characteristic polynomial of each of $X(\alpha)$, $X(\alpha)^{-1}$, $Y(\alpha)$, $Y(\alpha)^{-1}$ is
$$f(t) = (t-1)^2(t-1/\alpha^2)(t-\alpha^2)(t+\alpha)^2(t+1/\alpha)^2.$$
Assume that $||\alpha||>1$.  Then $f_1(t)=t-\alpha^2$ and $f_2(t)=f(t)/f_1(t)$.  An easy calculation shows the following.

\begin{enumerate}
\item $A(X(\alpha)) \leftrightarrow \text{span}\{-(\alpha+1)A_{31}+A_{32}\}$
\item $A(X(\alpha)^{-1}) \leftrightarrow \text{span}\{A_{12}-\frac{1+\alpha}{\alpha}A_{13}\}$
\item $A(Y(\alpha)) \leftrightarrow \text{span}\{-\alpha^2A_{21}+A_{23}-\alpha^2A_{31}-A_{32}-A_{22}\}$
\item $A(Y(\alpha)^{-1}) \leftrightarrow \text{span}\{-\frac{1}{\alpha}A_{12}-A_{13}+\alpha A_{21}-\alpha A_{23}+A_{11}\}$
\item $A'(X(\alpha)) \leftrightarrow \text{span}\{A_{23},A_{32},A_{22},A_{12},A_{13},A_{11},\frac{\alpha+1}{\alpha}A_{21}+A_{31}\}$
\item $A'(X(\alpha)^{-1})\leftrightarrow\text{span}\{A_{32},A_{22},A_{12},A_{11},A_{31},A_{21},A_{13}+(\alpha+1)A_{23}\}$
\item $A'(Y(\alpha))\leftrightarrow\text{span}\{-A_{31}+A_{22},A_{23},2A_{31}+A_{11},A_{12}+\alpha A_{31},A_{13},A_{21}-\frac{1}{\alpha}A_{31},
          -\alpha A_{31}+A_{32}\}$
\item $A'(Y(\alpha)^{-1})\leftrightarrow\text{span}\{A_{21},A_{12}-\alpha^2A_{23},-2A_{23}+A_{22},A_{23}+A_{32},A_{31},A_{13}-\alpha^2A_{23},A_{23}+A_{11}\}$.
\end{enumerate}
It is easy to check that conditions (\ref{cond2}) and (\ref{cond3}) hold.
\end{proof}

\begin{remark} That the matrices $\beta_4(\sigma_3\sigma_1^{-1})^m$ and $\beta_4(\sigma_2\sigma_3\sigma_1^{-1}\sigma_2^{-1})^m$
generate a free group was proved by S.~Moran \cite{moran} using the same technique over the field
$\zc(t)$.  It was our hope that passing to the matrices $X(\alpha),Y(\alpha)$ would allow us to take
$M=1$.  This is not the case however.  Indeed, denote by $v_1$ the basis vector of $A(X(\alpha))$ and by
$v_2$ the basis vector of $A(Y(\alpha))$.  Then it is easy to see that
$$||Yv_1-v_2||\ge 1 \qquad \text{and}\qquad ||Y^2v_1-v_2||\ge 1$$
so that no neighborhood of $A(X(\alpha))$ in ${\mathbb P}_V$ can be taken into a small neighborhood of
$A(Y(\alpha))$.  Proposition 1.1 of \cite{tits} therefore does not apply to $\langle X(\alpha),Y(\alpha)\rangle$.
\end{remark}

\subsection{Faithfulness of $\mu_n(\alpha,\beta)$}  The map $\kappa_n:B_n \ra GL_{\binom{n}{2}}(\zz[t^{\pm 1},q^{\pm 1}])$ is
faithful for all $n$.  In particular, if $\alpha$ and $\beta$ are algebraically independent complex numbers, then
the map induced by the homomorphism $t\mapsto \alpha$, $q\mapsto \beta$ yields a faithful representation
$$B_n\lra GL_{\binom{n}{2}}(\zz[\alpha^{\pm 1},\beta^{\pm 1}]).$$
Thus, $B_n\cap L(\alpha,\beta) = \{I\}$ in this case.

Recall that $\mu_n(\alpha,\beta)$ is the composition of $\kappa_n$ with the map 
$$g_n(\alpha,\beta):GL_{\binom{n}{2}}(\zz[t^{\pm 1},q^{\pm 1}])\lra \text{Aut}(L).$$
We note the following.

\begin{lemma} Denote by $Z$ the center of $GL_{\binom{n}{2}}(\zz[t^{\pm 1},q^{\pm 1}])$.
Then the kernel of $g_n(\alpha,\beta)$ is the subgroup $Z\cdot L(\alpha,\beta)$.
\end{lemma}

\begin{proof} It is clear that the kernel contains $Z\cdot L(\alpha,\beta)$.  For the reverse
inclusion, note that any $X\in GL_{\binom{n}{2}}(\zz[t^{\pm 1},q^{\pm 1}])$ can be written as
$X=X_0X_1$ where $X_0\in GL_{\binom{n}{2}}(\zc)$ and $X_1\in L(\alpha,\beta)$.  The action of $X$
on $L$ is then given by conjugation by $X_0$.  By considering the action on the basis of $L$, if $X_0$
acts trivially, we see that $X_0$ must be a diagonal matrix with all entries equal; that is, $X_0\in Z$.
\end{proof}

As $\kappa_n$ is faithful, we may identify $B_n$ with its image in $GL_{\binom{n}{2}}(\zz[t^{\pm 1},q^{\pm 1}])$.
Then if $\alpha$ and $\beta$ are algebraically independent, we see that the kernel of $\mu_n(\alpha,\beta)$
is $B_n\cap Z$.  We have thus proved the following result.

\begin{theorem}\label{mufaith} If $\alpha$ and $\beta$ are algebraically independent, then the kernel
of $\mu_n(\alpha,\beta)$ is precisely the center of $B_n$. \hfill $\qed$
\end{theorem}

\end{document}